\documentclass[submission]{FPSAC2020}

\usepackage[utf8x]{inputenc}
\usepackage[english]{babel}
\usepackage{amsmath,amsfonts,amssymb,amsthm}
\usepackage[T1]{fontenc}
\usepackage{hyperref}
\usepackage{cite}
\usepackage{dsfont}
\usepackage{paralist}
\usepackage{multicol}
\usepackage{stmaryrd}
\usepackage{subfig}
\usepackage{multirow}
\usepackage{tikz}
\usetikzlibrary{fit}
\usetikzlibrary{arrows.meta}

\usepackage{xcolor}
\definecolor{ColBlack}{RGB}{0,0,0} 
\definecolor{ColWhite}{RGB}{255,255,255} 
\definecolor{ColAA}{HTML}{520db1} 
\definecolor{ColAB}{HTML}{1a34c0} 
\definecolor{ColAC}{HTML}{3851db} 
\definecolor{ColBA}{HTML}{a80b3a} 
\definecolor{ColBB}{HTML}{a80b27} 
\definecolor{ColBC}{HTML}{b10d0d} 

\newtheorem{Theorem}{Theorem}[subsection]
\newtheorem{Proposition}[Theorem]{Proposition}

\numberwithin{equation}{section}

\renewcommand{\leq}{\leqslant}
\renewcommand{\geq}{\geqslant}

\title[Three Fuss-Catalan posets]{Three interacting families of Fuss-Catalan posets}
\author{
    Camille Combe%
    \thanks{\href{mailto:combe@math.unistra.fr}
        {\tt combe@math.unistra.fr}}%
    \addressmark{1}
    \and
    Samuele Giraudo%
    \thanks{\href{mailto:samuele.giraudo@u-pem.fr}
        {\tt samuele.giraudo@u-pem.fr}}%
    \addressmark{2}}
\address{
    \addressmark{1} Univ. Strasbourg, IRMA (UMR 7501), CNRS, France.
    \\
    \addressmark{2} LIGM, Univ. Gustave Eiffel, CNRS, ESIEE Paris, F-$77454$
    Marne-la-Vallée, France.}
\abstract{Three families of posets depending on a nonnegative integer parameter $m$ are
introduced. The underlying sets of these posets are enumerated by the $m$-Fuss Catalan
numbers. Among these, one is a generalization of Stanley lattices and another one is a
generalization of Tamari lattices. The three families of posets are related: they fit into a
chain for the order extension relation and they share some properties. Two associative
algebras are constructed as quotients of generalizations of the Malvenuto-Reutenauer
algebra. Their products describe intervals of our analogues of Stanley lattices and Tamari
lattices. In particular, one is a generalization of the Loday-Ronco algebra.}
\received{\today}

\tikzstyle{Centering}=[{baseline={([yshift=-0.5ex]current
    bounding box.center)}}]

\tikzstyle{MarkAA}=[draw=ColAA!80,fill=ColAA!8]
\tikzstyle{MarkAB}=[draw=ColAB!80,fill=ColAB!8]
\tikzstyle{MarkAC}=[draw=ColAC!80,fill=ColAC!8]
\tikzstyle{MarkBA}=[draw=ColBA!80,fill=ColBA!8]
\tikzstyle{MarkBB}=[draw=ColBB!80,fill=ColBB!8]
\tikzstyle{MarkBC}=[draw=ColBC!80,fill=ColBC!8]

\tikzstyle{Node}=[circle,MarkAA,inner sep=1pt,minimum size=2mm,thick,font=\scriptsize]
\tikzstyle{Edge}=[draw=ColBB!80,cap=round,thick,rounded corners=2.5pt]
\tikzstyle{Leaf}=[rectangle,MarkBC,inner sep=0pt,minimum size=1mm,thick]
\tikzstyle{NodeST}=[font=\footnotesize]
\tikzstyle{EdgeLabel}=[midway,inner sep=1pt,fill=ColWhite!0,font=\scriptsize]
\tikzstyle{LeafLabel}=[font=\scriptsize,node distance=2mm]
\tikzstyle{Subtree}=[regular polygon,regular polygon sides=3,MarkA,thick,minimum size=5mm,
font=\scriptsize]

\tikzstyle{PathNode}=[circle,MarkA,thick,inner sep=0pt,minimum size=2mm]
\tikzstyle{PathStep}=[color=Col1!60,thick]

\tikzstyle{Injection}=[ColBlack!100,draw,{>[scale=1.5,length=4,width=5]}-{>[scale=1.5,
length=4,width=5]}]
\tikzstyle{Surjection}=[ColBlack!100,draw,-{>[scale=1.5,length=4,width=5]>[scale=1.5,
length=4,width=5]}]
\tikzstyle{Map}=[ColBlack!100,draw,-{>[scale=1.5,length=4,width=5]}]

\tikzstyle{LineGrid}=[very thin,dashed,draw=ColAC!70]
\tikzstyle{Grid}=[LineGrid]

\tikzstyle{NodeGraph}=[circle,MarkAB,inner sep=1pt,minimum size=1.5mm,thick]
\tikzstyle{NodeLabeldGraph}=[font=\tiny,node distance=3mm]
\tikzstyle{EdgeGraph}=[ColBB!70,cap=round,very thick]
\tikzstyle{EdgeLabel}=[midway,inner sep=1pt,fill=ColWhite!0,font=\tiny]
\tikzstyle{FaceXY}=[fill=ColAA,opacity=.1]
\tikzstyle{FaceXZ}=[fill=ColBA,opacity=.2]
\tikzstyle{FaceYZ}=[fill=ColBC,opacity=.2]

\newcommand{\ColAA}[1]{\textcolor{ColAA}{#1}}
\newcommand{\ColAB}[1]{\textcolor{ColAB}{#1}}

\newcommand{\Hide}[1]{\ColAA{\tt HIDEN}}
\newcommand{\Def}[1]{\ColAB{\em #1}}
\newcommand{\Par}[1]{\left(#1\right)}
\newcommand{\Bra}[1]{\left\{#1\right\}}

\newcommand{\Han}[1]{\left[#1\right]}

\newcommand{\N}{\mathbb{N}}
\newcommand{\Z}{\mathbb{Z}}

\newcommand{\K}{\mathbb{K}}

\newcommand{\SetCliff}{\mathsf{Cl}}
\newcommand{\Weight}{\omega}
\newcommand{\Leq}{\preccurlyeq}
\newcommand{\PosetP}{\mathcal{P}}
\newcommand{\OutputWings}{{\mathcal{O}}}
\newcommand{\InputWings}{{\mathcal{I}}}
\newcommand{\Butterflies}{\mathcal{B}}
\newcommand{\Covered}{\lessdot}
\newcommand{\ElevationMap}{\mathbf{e}}
\newcommand{\ElevationImage}{\mathcal{E}}
\newcommand{\IncrMap}{{\Uparrow\!}}
\newcommand{\DecrMap}{{\Downarrow\!}}
\newcommand{\JJoin}{\vee}
\newcommand{\Meet}{\wedge}
\newcommand{\MaxLastLetter}{\mathrm{m}}
\newcommand{\DerivationOnWord}{\mathrm{d}}
\newcommand{\DerivationOnSet}{\mathrm{D}}

\newcommand{\Angle}[1]{\left\langle#1\right\rangle}
\newcommand{\SetAvalanche}{\mathsf{Av}}
\newcommand{\SetHill}{\mathsf{Hi}}
\newcommand{\SetCanyon}{\mathsf{Ca}}
\newcommand{\FussCatalan}{\mathrm{cat}}
\newcommand{\Reduction}{\mathrm{r}}

\newcommand{\Product}{\cdot}
\newcommand{\SpaceCliff}{\mathbf{Cl}}
\newcommand{\BasisE}{\mathsf{E}}
\newcommand{\BasisF}{\mathsf{F}}
\newcommand{\BasisH}{\mathsf{H}}
\newcommand{\SetPrime}{\mathcal{P}}
\newcommand{\AlphabetVar}{\mathbb{A}}
\newcommand{\RelationSpace}{\mathcal{R}}
\newcommand{\LeqSuffix}{\leq_{\mathrm{s}}}
\newcommand{\FQSym}{\mathsf{FQSym}}
\newcommand{\PBT}{\mathsf{PBT}}
\newcommand{\NCSym}{\mathsf{NCSym}}
\newcommand{\SpaceV}{\mathcal{V}}
\newcommand{\SpaceHill}{\mathbf{Hi}}
\newcommand{\SpaceCanyon}{\mathbf{Ca}}
\newcommand{\SubFamilly}{\mathcal{S}}
\newcommand{\IndicatorFunction}{\mathds{1}}

\DeclareMathOperator{\Over}{
\begin{tikzpicture}[Centering,scale=.16]
    \draw(0,0)--(1,1.25);
\end{tikzpicture}}

\DeclareMathOperator{\Under}{
\begin{tikzpicture}[Centering,scale=.16]
    \draw(0,0)--(1,-1.25);
\end{tikzpicture}}

\begin{document}
\maketitle

\section*{Introduction}
The theory of combinatorial Hopf algebras takes a prominent place in algebraic
combinatorics. The Malvenuto-Reutenauer algebra $\FQSym$~\cite{MR95,DHT02} is a central
object in this theory. This structure is defined on the linear span of all permutations and
the product of two permutations has the notable property to form an interval of the right
weak order. Moreover, $\FQSym$ admits a lot of substructures, like the Loday-Ronco algebra
of binary trees $\PBT$~\cite{LR98,HNT05} and the algebra of noncommutative symmetric
functions $\NCSym$~\cite{GKLLRT94}. Each of these structures bring out in a beautiful and
somewhat unexpected way the combinatorics of some partial orders, respectively the Tamari
order~\cite{Tam62} and the Boolean lattice, playing the same role as the one played by the
right weak order for~$\FQSym$.

The starting point of this work is to start from a different poset on permutations and ask
to what extent analogues of $\FQSym$ and a similar hierarchy of algebras arise in this
context. We consider here the componentwise ordering $\Leq$ on Lehmer codes of
permutations~\cite{Leh60}.  A study of these posets $\SetCliff_{\mathbf 1}(n)$ appears
in~\cite{Den13}. Each poset $\SetCliff_{\mathbf 1}(n)$ is an order extension of the right
weak order of order $n$.  The Hasse diagrams of the right weak order of order $3$ and of
$\SetCliff_{\mathbf 1}(3)$ are respectively

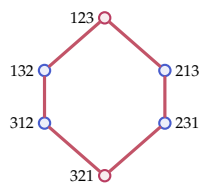
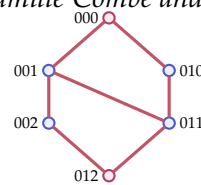
\begin{minipage}[t][0cm][b]{.4\textwidth}
\begin{equation}
    \scalebox{1}{
    \begin{tikzpicture}[Centering,xscale=.8,yscale=.7]
        \node[NodeGraph,MarkBA](000)[]at(0,0){};
        \node[NodeGraph](001)[]at(-1,-1){};
        \node[NodeGraph](002)[]at(-1,-2){};
        \node[NodeGraph](010)[]at(1,-1){};
        \node[NodeGraph](011)[]at(1,-2){};
        \node[NodeGraph,MarkBA](012)[]at(0,-3){};
        \node[NodeLabeldGraph,left of=000]{$123$};
        \node[NodeLabeldGraph,left of=001]{$132$};
        \node[NodeLabeldGraph,left of=002]{$312$};
        \node[NodeLabeldGraph,right of=010]{$213$};
        \node[NodeLabeldGraph,right of=011]{$231$};
        \node[NodeLabeldGraph,left of=012]{$321$};
        \draw[EdgeGraph](000)--(001);
        \draw[EdgeGraph](000)--(010);
        \draw[EdgeGraph](001)--(002);
        \draw[EdgeGraph](010)--(011);
        \draw[EdgeGraph](002)--(012);
        \draw[EdgeGraph](011)--(012);
    \end{tikzpicture}}
\end{equation}
\end{minipage}
\qquad and
\begin{minipage}[t][0cm][b]{.4\textwidth}
\begin{equation}
    \scalebox{1}{
    \begin{tikzpicture}[Centering,xscale=.8,yscale=.7]
        \node[NodeGraph,MarkBA](000)[]at(0,0){};
        \node[NodeGraph](001)[]at(-1,-1){};
        \node[NodeGraph](002)[]at(-1,-2){};
        \node[NodeGraph](010)[]at(1,-1){};
        \node[NodeGraph](011)[]at(1,-2){};
        \node[NodeGraph,MarkBA](012)[]at(0,-3){};
        \node[NodeLabeldGraph,left of=000]{$000$};
        \node[NodeLabeldGraph,left of=001]{$001$};
        \node[NodeLabeldGraph,left of=002]{$002$};
        \node[NodeLabeldGraph,right of=010]{$010$};
        \node[NodeLabeldGraph,right of=011]{$011$};
        \node[NodeLabeldGraph,left of=012]{$012$};
        \draw[EdgeGraph](000)--(001);
        \draw[EdgeGraph](000)--(010);
        \draw[EdgeGraph](001)--(002);
        \draw[EdgeGraph](001)--(011);
        \draw[EdgeGraph](010)--(011);
        \draw[EdgeGraph](002)--(012);
        \draw[EdgeGraph](011)--(012);
    \end{tikzpicture}}.
\end{equation}
\end{minipage}

\noindent In this work, we consider a more general version of Lehmer codes, called
$\delta$-cliffs, leading to distributive lattices $\SetCliff_\delta$. Here $\delta$ is a
parameter which is a map $\N \setminus \{0\} \to \N$. The linear spans $\SpaceCliff_\delta$
of these sets are endowed with a very natural product related to the intervals of
$\SetCliff_\delta$, forming associative algebras when $\delta$ satisfies some property that
is called being valley-free. The following table inventories some properties of
$\SpaceCliff_\delta$ implied by properties of~$\delta$:
\begin{center} \footnotesize
    \begin{tabular}{|c|c|} \hline
        {\bf Properties of $\delta$} & {\bf Properties of $\SpaceCliff_\delta$}
            \\ \hline \hline
        \multirow{2}{*}{None} & Unital graded magmatic algebra \\
        & Products on the $\BasisF$-basis are intervals in
            $\delta$-cliff posets \\ \hline
        Valley-free & Associative algebra \\ \hline
        Valley-free and $1$-dominated & Finite presentation \\ \hline
        Weakly increasing & Free as unital associative algebra \\ \hline
    \end{tabular}
\end{center}
The particular algebra $\SpaceCliff_{\mathbf 1}$ is in fact isomorphic to
$\FQSym$, so that $\SpaceCliff_\delta$ is a generalization of this latter.

In the same way as the Tamari order can be defined by restricting the right weak order to
some permutations, one builds three subposets of $\SetCliff_\delta$ by restricting
$\Leq$ to particular $\delta$-cliffs. This leads to three families $\SetAvalanche_\delta$,
$\SetHill_\delta$, and $\SetCanyon_\delta$ of posets. When $\delta$ is a particular map
${\mathbf m}$ depending on an integer $m \geq 0$, the underlying sets of all these posets of
order $n \geq 0$ are enumerated by the $n$-th $m$-Fuss-Catalan number~\cite{DM47}

\begin{footnotesize}
\begin{equation}
    \FussCatalan_m(n) := \frac{1}{mn + 1} \binom{mn + n}{n}.
\end{equation}
\end{footnotesize}
These posets have some close interactions: $\SetHill_\delta$ is an order extension of
$\SetCanyon_\delta$, which is itself an order extension of $\SetAvalanche_\delta$. Besides,
$\SetHill_{\mathbf 1}$ (resp. $\SetCanyon_{\mathbf 1}$) is the Stanley lattice~\cite{Knu04}
(resp. the Tamari lattice), so that $\SetHill_{\mathbf m}$ (resp. $\SetCanyon_{\mathbf m}$),
$m \geq 0$, are new generalizations of Stanley lattices (resp. Tamari lattices
---see~\cite{BPR12} for the classical one).  From these posets $\SetHill_{\mathbf m}$ and
$\SetCanyon_{\mathbf m}$, one defines respectively two quotients $\SpaceHill_m$ and
$\SpaceCanyon_m$ of $\SpaceCliff_{\mathbf m}$. The algebra $\SpaceCanyon_1$ is isomorphic to
$\PBT$, and the other ones $\SpaceCanyon_m$, $m \geq 2$, are not free as associative
algebras.

This paper is organized as follows. Section~\ref{sec:cliff_posets} is intended to introduce
$\delta$-cliffs, the posets $\SetCliff_\delta$, and sufficient conditions on subposets
$\SubFamilly$ of $\SetCliff_\delta$ to satisfy some order theoretic properties such as
EL-shellability, lattice property, and constructibility by interval doubling. In
Section~\ref{sec:fuss_catalan_posets}, we study the posets $\SetAvalanche_\delta$,
$\SetHill_\delta$, and $\SetCanyon_\delta$, and their interactions. Finally,
Section~\ref{sec:algebras} presents a study of the algebras $\SpaceCliff_\delta$,
$\SpaceHill_m$, and $\SpaceCanyon_m$.

\vspace{-1em}

\paragraph{General notations and conventions.}
For any integers $i$ and $j$, $[i, j]$ denotes the set $\{i, i + 1, \dots, j\}$ and $[i]$
denotes the set $[1, i]$.  Graded sets are sets decomposing as a disjoint union $\SubFamilly
= \bigsqcup_{n \geq 0} \SubFamilly(n)$. The degree of $x \in \SubFamilly$ is the unique $n
\geq 0$ such that $x \in \SubFamilly(n)$. A graded subset of $\SubFamilly$ is a graded set
$\SubFamilly'$ such that for all $n \geq 0$, $\SubFamilly'(n) \subseteq \SubFamilly(n)$. If
$P$ is a statement, we denote by $\IndicatorFunction_P$ the indicator function (equals to
$1$ if $P$ holds and $0$ otherwise).

\section{Cliff posets and general properties} \label{sec:cliff_posets}

\subsection{Range maps and cliffs posets}
A \Def{range map} is a map $\delta : \N \setminus \{0\} \to \N$.  We shall specify range
maps as infinite words $\delta = \delta(1) \delta(2) \dots$. To this purpose, for any $a \in
\N$, we denote by $a^\omega$ the infinite word having all its letters equal to $a$.
Moreover, for any $m \geq 0$, we denote by $\mathbf{m}$ the range map satisfying $\mathbf{m}
:= 0 \, m \, (2m) \, (3m) \dots$.  We say that $\delta$ is \Def{weakly increasing} (resp.
\Def{increasing}) if for all $i \geq 1$, $\delta(i) \leq \delta(i + 1)$ (resp. $\delta(i) <
\delta(i + 1)$); that $\delta$ is \Def{valley-free} (or \Def{unimodal}) if there are no $1
\leq i_1 < i_2 < i_3$ such that $\delta\Par{i_1} > \delta\Par{i_2} < \delta\Par{i_3}$; and
that $\delta$ is \Def{$j$-dominated} for a $j \geq 1$ if there is $k \geq 1$ such that for
all $k' \geq k$, $\delta(j) \geq \delta\Par{k'}$.

Given a range map $\delta$, a word $u$ of integers of length $n$ is a \Def{$\delta$-cliff}
if for any $i \in [n]$, $0 \leq u_i \leq \delta(i)$. The \Def{size} $|u|$ of a
$\delta$-cliff $u$ is its length as a word, and the \Def{weight} $\Weight(u)$ of $u$ is the
sum of its letters. The graded set of all $\delta$-cliffs where the degree of a
$\delta$-cliff is its size, is denoted by $\SetCliff_\delta$. The set of
$\mathbf{1}$-cliffs of size $n$ is in one-to-one correspondence with the set of
permutations of the same size.  A possible correspondence sends a permutation $\sigma$ to
its Lehmer code~\cite{Leh60}. Moreover, when $\delta$ is weakly increasing, there is a
one-to-one correspondence between $\SetCliff_\delta(n)$ and the set of $s$-decreasing
trees~\cite{CP19}, where $s$ is obtained from $\delta$ and $n$. Therefore, one can see
$\delta$-cliffs as generalizations of both permutations and decreasing trees.

Let $\Leq$ be the partial order relation on $\SetCliff_\delta$ defined by $u \Leq v$ for any
$u, v \in \SetCliff_\delta$ such that $|u| = |v|$ and $u_i \leq v_i$ for all $i \in [|u|]$.
For any $n \geq 0$, the poset $\Par{\SetCliff_\delta(n), \Leq}$ is the \Def{$\delta$-cliff
poset} of order $n$.  A study of the $\mathbf{1}$-cliff posets appears in~\cite{Den13}. Our
definition stated here depending on $\delta$ is therefore a generalization of these posets.
The structure of the $\delta$-cliff posets is very simple since each of these posets of
order $n$ is isomorphic to the Cartesian product $[0, \delta(1)] \times \dots \times
[0, \delta(n)]$, where $[k]$ is the total order on $k$ elements. It follows from this
observation that each $\delta$-cliff poset is a lattice admitting respectively $\Meet$ and
$\JJoin$ as meet and join operations, defined respectively as the componentwise minimum and
maximum.

\subsection{Subposets of cliff posets}
Despite their simplicity, the $\delta$-cliff posets contain subposets having a lot of
combinatorial and algebraic properties.  If $\SubFamilly$ is a graded subset of
$\SetCliff_\delta$, each $\SubFamilly(n)$, $n \geq 0$, is a subposet of
$\SetCliff_\delta(n)$ for the order relation $\Leq$.  We denote by $\Covered_\SubFamilly$
the covering relation of $\SubFamilly(n)$. We say that $\SubFamilly$ is \Def{straight} if
for any $n \geq 0$ and any $u, v \in \SubFamilly(n)$ such that $u \Covered_\SubFamilly v$,
then $u$ and $v$ differ by just one letter; that $\SubFamilly$ is \Def{coated} if for any
$n \geq 0$, any $u, v \in \SubFamilly(n)$ such that $u \Leq v$, and any $i \in [n - 1]$,
then $u_1 \dots u_i v_{i + 1} \dots v_n \in \SubFamilly$; that $\SubFamilly$ is \Def{closed
by prefix} if for any $n \geq 0$ and any $u \in \SubFamilly(n)$, all prefixes of $u$ belong
to $\SubFamilly$; and that $\SubFamilly$ is \Def{minimally extendable} (resp. \Def{maximally
extendable}) if $\epsilon \in \SubFamilly$ and for any $u \in \SubFamilly$,
$u0 \in \SubFamilly$ (resp. $u \, \delta(n + 1) \in \SubFamilly$), where $\epsilon$ is the
$\delta$-cliff of size $0$.  When $\SubFamilly$ is straight, let for any $n \geq 0$,
$\alpha_n := \# \Bra{i \in [n] : \delta(i) \ne 0}$.  In this case, the graded set
$\InputWings(\SubFamilly)$ (resp. $\OutputWings(\SubFamilly)$) of \Def{input-wings} (resp.
\Def{output-wings}) contains all $u \in \SubFamilly$ which cover (resp. are covered by)
exactly $\alpha_{|u|}$ elements. Moreover, the graded set $\Butterflies(\SubFamilly)$ of
\Def{butterflies} is the intersection $\InputWings(\SubFamilly) \cap
\OutputWings(\SubFamilly)$.  We present here general results about subposets of
$\delta$-cliff posets.

\vspace{-1em}
\paragraph{Shellability.}
The next result concerns a sufficient condition on $\SubFamilly(n)$, $n \geq 0$, to have the
property to be EL-shellable~\cite{BW96}.

\begin{Theorem} \label{thm:subposets_el_shellability}
    Let $\SubFamilly$ be a straight and coated graded subset of $\SetCliff_\delta$.  For any
    $n \geq 0$, the map
    \begin{math}
        \lambda_\SubFamilly : \Covered_\SubFamilly \to \Z^2
    \end{math}
    defined for any $(u, v) \in \Covered_\SubFamilly$ by
    \begin{math}
        \lambda_\SubFamilly(u, v) := \Par{-i, u_i}
    \end{math}
    where $i$ is the unique position such $u_i \ne v_i$, is an EL-labeling of
    $\SubFamilly(n)$.  Moreover, there is at most one weakly decreasing chain between any
    pair of elements of~$\SubFamilly(n)$.
\end{Theorem}

By Theorem~\ref{thm:subposets_el_shellability}, when moreover $\SubFamilly$ is such that for
any $n \geq 0$, $\SubFamilly(n)$ has a least and greatest element, $\SubFamilly(n)$ is
EL-shellable.

\vspace{-1em}
\paragraph{Meet and join operations, and lattices.}
The next result provides a sufficient condition on $\SubFamilly(n)$, $n \geq 0$, to be a
lattice and describes an algorithm to compute its meet and join operations.  Assume that
$\SubFamilly$ is minimally extendable.  For any $n \geq 0$, the
\Def{$\SubFamilly$-decrementation map} is the map
\begin{math}
    \DecrMap_\SubFamilly : \SetCliff_\delta(n) \to \SubFamilly(n)
\end{math}
defined recursively by $\DecrMap_\SubFamilly(\epsilon) := \epsilon$ and, for any $u a \in
\SetCliff_\delta(n)$ where $u \in \SetCliff_\delta$ and $a \in \N$, by
\begin{math}
    \DecrMap_\SubFamilly(u a) := \DecrMap_\SubFamilly(u) \, b
\end{math}
where
\begin{math}
    b := \max \Bra{b \leq a : \DecrMap_\SubFamilly(u) \, b \in \SubFamilly}.
\end{math}
Observe that the fact that $\SubFamilly$ is minimally extendable ensures that
$\DecrMap_\SubFamilly$ is a well-defined map. Let also, for any $n \geq 0$ and $u, v \in
\SubFamilly(n)$,
\begin{math}
    u \Meet_\SubFamilly v := \DecrMap_\SubFamilly(u \Meet v).
\end{math}
When $\SubFamilly$ is maximally extendable, we define dually the
\Def{$\SubFamilly$-incrementation map} $\IncrMap_\SubFamilly$ and the operation
$\JJoin_\SubFamilly$.

\begin{Theorem} \label{thm:decr_incr_meet_join}
    Let $\SubFamilly$ be a closed by prefix graded subset of $\SetCliff_\delta$. If
    $\SubFamilly$ is minimally (resp. maximally) extendable, then the operation
    $\Meet_\SubFamilly$ (resp.  $\JJoin_\SubFamilly$) is, for any $n \geq 0$, the meet
    (resp. join) operation of the poset $\SubFamilly(n)$.
\end{Theorem}

\vspace{-1em}
\paragraph{Constructibility by interval doubling.}
The next result brings a sufficient condition for $\SubFamilly(n)$, $n \geq 0$, to be
constructible by interval doubling~\cite{Day92} and describes explicitly the sequence of the
involved interval contractions (the inverse operation of interval doubling).

Let $\PosetP$ be a nonempty subposet of $\SetCliff_\delta(n)$ for a given fixed size $n \geq
0$.  Let $\MaxLastLetter(\PosetP) := \max \Bra{u_n : u \in \PosetP}$. When
$\MaxLastLetter(\PosetP) \geq 1$, let $\DerivationOnWord_\PosetP : \PosetP \to
\SetCliff_\delta(n)$ be the map sending any $u \in \PosetP$ to the word obtained by
decrementing by $1$ the last letter $u_n$ of $u$ if $u_n = \MaxLastLetter(\PosetP)$, and to
$u$ otherwise. The \Def{derivation} of $\PosetP$ is the set $\DerivationOnSet(\PosetP) :=
\Bra{\DerivationOnWord_\PosetP(u) : u \in \PosetP}$.  Observe that
$\DerivationOnSet(\PosetP)$ is a subposet of $\SetCliff_\delta(n)$ but not necessarily a
subposet of $\PosetP$. Moreover, for any $a \in [0, \delta(n)]$, let $\PosetP_a := \Bra{u
\in \PosetP : u_n = a}$ and
\begin{math}
    \PosetP_a^\MaxLastLetter
    := \Bra{u_1 \dots u_{n - 1} \; \MaxLastLetter(\PosetP)
    : u \in \PosetP_a}.
\end{math}
Observe that $\PosetP_a$ is a subposet of $\PosetP$ while $\PosetP_a^\MaxLastLetter$ may
contain $\delta$-cliffs that do not belong to $\PosetP$.  The subposet $\PosetP$ is
\Def{nested} if for any $a \in \Han{0, \MaxLastLetter(\PosetP)}$, the $\delta$-cliff $0^{n -
1} a$ belongs to $\PosetP$ and $\PosetP_a^\MaxLastLetter$ is both a subset and an interval
of~$\PosetP$.  This definition still holds when $\MaxLastLetter(\PosetP) = 0$. Observe that
any $\delta$-cliff $0^{n - 1} a$, $a \geq 1$, of $\PosetP$ covers exactly the single element
$0^{n - 1} (a - 1)$ of $\PosetP$. Therefore, when $\PosetP$ is a lattice, these
$\delta$-cliffs are join-irreducible.  We say by extension that a graded subset
$\SubFamilly$ of $\SetCliff_\delta$ is \Def{nested} if for all $n \geq 0$, the posets
$\SubFamilly(n)$ are nested.

\begin{Theorem} \label{thm:constructible_by_interval_doubling_subfamilly}
    Let $\SubFamilly$ be a nested and closed by prefix graded subset of $\SetCliff_\delta$
    where $\delta(1) = 0$.  For any $n \geq 1$, $\SubFamilly(n)$ is constructible by
    interval doubling. Moreover,
    \begin{small}
    \begin{equation} \begin{split}
        \SubFamilly(n) & \to \DerivationOnSet(\SubFamilly(n))
        \to \dots \to \DerivationOnSet^{\MaxLastLetter(\SubFamilly(n))}(\SubFamilly(n))
        \simeq \SubFamilly(n - 1)
        \\
        & \to \DerivationOnSet(\SubFamilly(n - 1)) \to \dots
        \to \DerivationOnSet^{\MaxLastLetter(\SubFamilly(n - 1))}(\SubFamilly(n - 1))
        \simeq \SubFamilly(n - 2) \\
        & \to \dots \to \SubFamilly(0) \simeq \{\epsilon\}
    \end{split} \end{equation}
    \end{small}
    is a sequence of interval contractions from $\SubFamilly(n)$ to the trivial
    lattice~$\{\epsilon\}$.
\end{Theorem}

\paragraph{Elevation map.}
We introduce here a combinatorial tool intervening in the study of the three Fuss-Catalan
posets introduced in the next section.

When $\SubFamilly$ is closed by prefix, for any $n \geq 0$, the
\Def{$\SubFamilly$-elevation map} is the map
\begin{math}
    \ElevationMap_\SubFamilly : \SubFamilly(n) \to \SetCliff_\delta(n)
\end{math}
defined, for any $u \in \SubFamilly(n)$ by setting
\begin{equation}
    \ElevationMap_\SubFamilly(u)_i
    :=
    \# \Bra{a \in \N  : a < u_i \mbox{ and } u_1 \dots u_{i - 1} a \in \SubFamilly}
\end{equation}
for any $i \in [n]$. From an intuitive point of view, the value of the $i$-th letter of
$\ElevationMap_\SubFamilly(u)$ is the number of words of $\SubFamilly$ obtained by
considering the prefix of $u$ ending at the letter $u_i$ and by replacing this letter by a
smaller one. Observe in particular that $\ElevationMap_{\SetCliff_\delta}$ is the
identity map. Let $\ElevationImage_\SubFamilly$ be the graded set wherein for any
$n \geq 0$, $\ElevationImage_\SubFamilly(n)$ is the image of $\SubFamilly(n)$ by the
$\SubFamilly$-elevation map.  We call this set the \Def{$\SubFamilly$-elevation image}.
Observe that $\ElevationImage_\SubFamilly$ is a graded subset of $\SetCliff_\delta$ and that
for any $u \in \SubFamilly$, $\ElevationMap_\SubFamilly(u) \Leq u$.

\begin{Proposition} \label{prop:elevation_map}
    Let $\SubFamilly$ be a closed by prefix graded subset of $\SetCliff_\delta$. Then,
    \begin{inparaenum}[(i)]
        \item \label{item:elevation_map_1} the $\SubFamilly$-elevation map is injective;
        \item the $\SubFamilly$-elevation image is closed by prefix;
        \item \label{item:elevation_map_3}
        if for any $u, v \in \SubFamilly$ such that $u \Leq v$, $v a \in \SubFamilly$
        implies $u a \in \SubFamilly$ where $a \in \N$, then for any $n \geq 0$, the map
        $\ElevationMap_\SubFamilly^{-1}$ is a poset morphism from
        $\ElevationImage_\SubFamilly(n)$ to~$\SubFamilly(n)$.
    \end{inparaenum}
\end{Proposition}
Point~\eqref{item:elevation_map_3} of Proposition~\ref{prop:elevation_map} says that when
$\SubFamilly$ satisfies the given prerequisites, every poset $\SubFamilly(n)$ is an order
extension of $\ElevationImage_\SubFamilly(n)$.

\section{Three Fuss-Catalan subposets} \label{sec:fuss_catalan_posets}

\subsection{Avalanche posets}
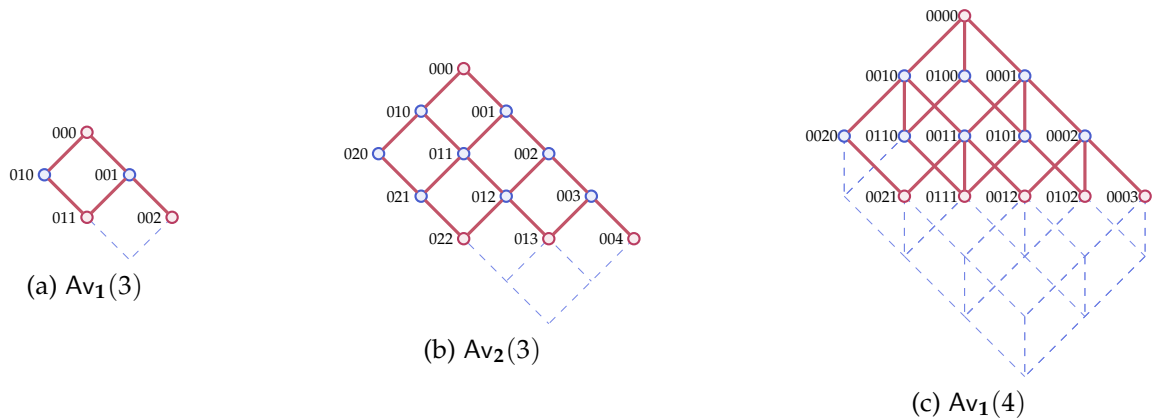
\begin{figure}[ht]
    \centering
    \subfloat[][$\SetAvalanche_\mathbf{1}(3)$]{
    \centering
    \scalebox{1}{
    \begin{tikzpicture}[Centering,xscale=.8,yscale=.8,rotate=-135]
        \draw[Grid](0,0)grid(1,2);
        \node[NodeGraph,MarkBA](000)[]at(0,0){};
        \node[NodeGraph](001)[]at(0,1){};
        \node[NodeGraph,MarkBA](002)[]at(0,2){};
        \node[NodeGraph](010)[]at(,0){};
        \node[NodeGraph,MarkBA](011)[]at(1,1){};
        \node[NodeLabeldGraph,left of=000]{$000$};
        \node[NodeLabeldGraph,left of=001]{$001$};
        \node[NodeLabeldGraph,left of=002]{$002$};
        \node[NodeLabeldGraph,left of=010]{$010$};
        \node[NodeLabeldGraph,left of=011]{$011$};
        \draw[EdgeGraph](000)--(001);
        \draw[EdgeGraph](000)--(010);
        \draw[EdgeGraph](001)--(002);
        \draw[EdgeGraph](001)--(011);
        \draw[EdgeGraph](010)--(011);
    \end{tikzpicture}}
    \label{subfig:avalanche_poset_1_3}}
    \qquad \qquad
    \subfloat[][$\SetAvalanche_\mathbf{2}(3)$]{
    \centering
    \scalebox{1}{
    \begin{tikzpicture}[Centering,xscale=.8,yscale=.8,rotate=-135]
        \draw[Grid](0,0)grid(2,4);
        \node[NodeGraph,MarkBA](000)[]at(0,0){};
        \node[NodeGraph](001)[]at(0,1){};
        \node[NodeGraph](002)[]at(0,2){};
        \node[NodeGraph](003)[]at(0,3){};
        \node[NodeGraph,MarkBA](004)[]at(0,4){};
        \node[NodeGraph](010)[]at(1,0){};
        \node[NodeGraph](011)[]at(1,1){};
        \node[NodeGraph](012)[]at(1,2){};
        \node[NodeGraph,MarkBA](013)[]at(1,3){};
        \node[NodeGraph](020)[]at(2,0){};
        \node[NodeGraph](021)[]at(2,1){};
        \node[NodeGraph,MarkBA](022)[]at(2,2){};
        \node[NodeLabeldGraph,left of=000]{$000$};
        \node[NodeLabeldGraph,left of=001]{$001$};
        \node[NodeLabeldGraph,left of=002]{$002$};
        \node[NodeLabeldGraph,left of=003]{$003$};
        \node[NodeLabeldGraph,left of=004]{$004$};
        \node[NodeLabeldGraph,left of=010]{$010$};
        \node[NodeLabeldGraph,left of=011]{$011$};
        \node[NodeLabeldGraph,left of=012]{$012$};
        \node[NodeLabeldGraph,left of=013]{$013$};
        \node[NodeLabeldGraph,left of=020]{$020$};
        \node[NodeLabeldGraph,left of=021]{$021$};
        \node[NodeLabeldGraph,left of=022]{$022$};
        \draw[EdgeGraph](000)--(001);
        \draw[EdgeGraph](000)--(010);
        \draw[EdgeGraph](001)--(002);
        \draw[EdgeGraph](001)--(011);
        \draw[EdgeGraph](002)--(003);
        \draw[EdgeGraph](002)--(012);
        \draw[EdgeGraph](003)--(004);
        \draw[EdgeGraph](003)--(013);
        \draw[EdgeGraph](010)--(011);
        \draw[EdgeGraph](010)--(020);
        \draw[EdgeGraph](011)--(012);
        \draw[EdgeGraph](011)--(021);
        \draw[EdgeGraph](012)--(013);
        \draw[EdgeGraph](012)--(022);
        \draw[EdgeGraph](020)--(021);
        \draw[EdgeGraph](021)--(022);
        \end{tikzpicture}}
    \label{subfig:avalanche_poset_2_3}}
    \qquad \qquad
    \subfloat[][$\SetAvalanche_\mathbf{1}(4)$]{
    \centering
    \scalebox{1}{
    \begin{tikzpicture}[Centering,xscale=.8,yscale=.8,
        x={(0,-1cm)},
        y={(-1cm,-1cm)},
        z={(1cm,-1cm)}]
        \foreach \x in {0,...,1} {
            \foreach \y in {0,...,2} {
                \foreach \z in {0,..., 3} {
                    \draw[LineGrid](\x,0,\z)--(\x,2,\z);
                    \draw[LineGrid](0,\y,\z)--(1,\y,\z);
                }
                \foreach \z in {0,..., 2} {
                    \draw[LineGrid](\x,\y,\z)--(\x,\y,\z+1);
                }
            }
        }
        \node[NodeGraph,MarkBA](0000)[]at(0,0,0){};
        \node[NodeGraph](0001)[]at(0,0,1){};
        \node[NodeGraph](0002)[]at(0,0,2){};
        \node[NodeGraph,MarkBA](0003)[]at(0,0,3){};
        \node[NodeGraph](0010)[]at(0,1,0){};
        \node[NodeGraph](0011)[]at(0,1,1){};
        \node[NodeGraph,MarkBA](0012)[]at(0,1,2){};
        \node[NodeGraph](0020)[]at(0,2,0){};
        \node[NodeGraph,MarkBA](0021)[]at(0,2,1){};
        \node[NodeGraph](0100)[]at(1,0,0){};
        \node[NodeGraph](0101)[]at(1,0,1){};
        \node[NodeGraph,MarkBA](0102)[]at(1,0,2){};
        \node[NodeGraph](0110)[]at(1,1,0){};
        \node[NodeGraph,MarkBA](0111)[]at(1,1,1){};
        \node[NodeLabeldGraph,left of=0000]{$0000$};
        \node[NodeLabeldGraph,left of=0001]{$0001$};
        \node[NodeLabeldGraph,left of=0002]{$0002$};
        \node[NodeLabeldGraph,left of=0003]{$0003$};
        \node[NodeLabeldGraph,left of=0010]{$0010$};
        \node[NodeLabeldGraph,left of=0011]{$0011$};
        \node[NodeLabeldGraph,left of=0012]{$0012$};
        \node[NodeLabeldGraph,left of=0020]{$0020$};
        \node[NodeLabeldGraph,left of=0021]{$0021$};
        \node[NodeLabeldGraph,left of=0100]{$0100$};
        \node[NodeLabeldGraph,left of=0101]{$0101$};
        \node[NodeLabeldGraph,left of=0102]{$0102$};
        \node[NodeLabeldGraph,left of=0110]{$0110$};
        \node[NodeLabeldGraph,left of=0111]{$0111$};
        \draw[EdgeGraph](0000)--(0001);
        \draw[EdgeGraph](0000)--(0010);
        \draw[EdgeGraph](0000)--(0100);
        \draw[EdgeGraph](0001)--(0002);
        \draw[EdgeGraph](0001)--(0011);
        \draw[EdgeGraph](0001)--(0101);
        \draw[EdgeGraph](0002)--(0003);
        \draw[EdgeGraph](0002)--(0012);
        \draw[EdgeGraph](0002)--(0102);
        \draw[EdgeGraph](0010)--(0011);
        \draw[EdgeGraph](0010)--(0020);
        \draw[EdgeGraph](0010)--(0110);
        \draw[EdgeGraph](0011)--(0012);
        \draw[EdgeGraph](0011)--(0021);
        \draw[EdgeGraph](0011)--(0111);
        \draw[EdgeGraph](0020)--(0021);
        \draw[EdgeGraph](0100)--(0101);
        \draw[EdgeGraph](0100)--(0110);
        \draw[EdgeGraph](0101)--(0102);
        \draw[EdgeGraph](0101)--(0111);
        \draw[EdgeGraph](0110)--(0111);
    \end{tikzpicture}}
    \label{subfig:avalanche_poset_1_4}}
    \caption{\footnotesize Hasse diagrams of some $\delta$-avalanche posets.}
    \label{fig:examples_avalanche_posets}
\end{figure}

Let $\SetAvalanche_\delta$ be the graded subset of $\SetCliff_\delta$ containing all
$\delta$-cliffs $u$ such that for all nonempty prefixes $u'$ of $u$, $\Weight\Par{u'} \leq
\delta\Par{\left|u'\right|}$. Any element of $\SetAvalanche_\delta$ is a
\Def{$\delta$-avalanche}.  Figure~\ref{fig:examples_avalanche_posets} shows some ${\mathbf
m}$-avalanche posets.
Notice that there are several maximal elements in $\SetAvalanche_{\mathbf m}(n)$.
Obviously, since by definition we have in particular, for all $u \in \SetAvalanche_{\mathbf
m}(n)$, $\Weight\Par{u} \leq \delta\Par{\left|u\right|}$, then $u$ is a maximal element of
$\SetAvalanche_{\mathbf m}(n)$ if and only if $\Weight\Par{u} = m (n - 1)$.

\begin{Proposition} \label{prop:properties_avalanche_objects}
    For any $m \geq 0$, $\SetAvalanche_{\mathbf m}$ is straight, coated, closed by prefix,
    and minimally extendable. Moreover, for any $n \geq 0$, $\SetAvalanche_{\mathbf m}(n)$
    is graded and a meet semi-sublattice of $\SetCliff_{\mathbf m}(n)$.
\end{Proposition}

\begin{Proposition} \label{prop:enumeration_avalanche}
    For any $m \geq 0$ and $n \geq 0$,
    \begin{math}
        \# \SetAvalanche_{\mathbf m}(n) = \FussCatalan_m(n).
    \end{math}
\end{Proposition}

Proposition~\ref{prop:enumeration_avalanche} is a consequence of the fact that there is a
bijection between ${\mathbf m}$-avalanches and $m$-Dyck paths, objects enumerated by
Fuss-Catalan numbers.

\subsection{Hill lattices} \label{sec:hill_lattices}
\begin{figure}[ht]
    \centering
    \subfloat[][$\SetHill_\mathbf{1}(3)$]{
    \centering
    \scalebox{1}{
    \begin{tikzpicture}[Centering,xscale=.8,yscale=.8,rotate=-135]
        \draw[Grid](0,0)grid(1,2);
        \node[NodeGraph,MarkBA](000)[]at(0,0){};
        \node[NodeGraph](001)[]at(0,1){};
        \node[NodeGraph](002)[]at(0,2){};
        \node[NodeGraph](011)[]at(1,1){};
        \node[NodeGraph,MarkBA](012)[]at(1,2){};
        \node[NodeLabeldGraph,left of=000]{$000$};
        \node[NodeLabeldGraph,left of=001]{$001$};
        \node[NodeLabeldGraph,left of=002]{$002$};
        \node[NodeLabeldGraph,left of=011]{$011$};
        \node[NodeLabeldGraph,left of=012]{$012$};
        \draw[EdgeGraph](000)--(001);
        \draw[EdgeGraph](001)--(002);
        \draw[EdgeGraph](001)--(011);
        \draw[EdgeGraph](002)--(012);
        \draw[EdgeGraph](011)--(012);
    \end{tikzpicture}}
    \label{subfig:hill_poset_1_3}}
    \qquad \qquad
    \subfloat[][$\SetHill_\mathbf{2}(3)$]{
    \centering
    \scalebox{1}{
    \begin{tikzpicture}[Centering,xscale=.8,yscale=.8,rotate=-135]
        \draw[Grid](0,0)grid(2,4);
        \node[NodeGraph,MarkBA](000)[]at(0,0){};
        \node[NodeGraph](001)[]at(0,1){};
        \node[NodeGraph](002)[]at(0,2){};
        \node[NodeGraph](003)[]at(0,3){};
        \node[NodeGraph](004)[]at(0,4){};
        \node[NodeGraph](011)[]at(1,1){};
        \node[NodeGraph](012)[]at(1,2){};
        \node[NodeGraph](013)[]at(1,3){};
        \node[NodeGraph](014)[]at(1,4){};
        \node[NodeGraph](022)[]at(2,2){};
        \node[NodeGraph](023)[]at(2,3){};
        \node[NodeGraph,MarkBA](024)[]at(2,4){};
        \node[NodeLabeldGraph,left of=000]{$000$};
        \node[NodeLabeldGraph,left of=001]{$001$};
        \node[NodeLabeldGraph,left of=002]{$002$};
        \node[NodeLabeldGraph,left of=003]{$003$};
        \node[NodeLabeldGraph,left of=004]{$004$};
        \node[NodeLabeldGraph,left of=011]{$011$};
        \node[NodeLabeldGraph,left of=012]{$012$};
        \node[NodeLabeldGraph,left of=013]{$013$};
        \node[NodeLabeldGraph,left of=014]{$014$};
        \node[NodeLabeldGraph,left of=022]{$022$};
        \node[NodeLabeldGraph,left of=023]{$023$};
        \node[NodeLabeldGraph,left of=024]{$024$};
        \draw[EdgeGraph](000)--(001);
        \draw[EdgeGraph](001)--(002);
        \draw[EdgeGraph](001)--(011);
        \draw[EdgeGraph](002)--(003);
        \draw[EdgeGraph](002)--(012);
        \draw[EdgeGraph](003)--(004);
        \draw[EdgeGraph](003)--(013);
        \draw[EdgeGraph](004)--(014);
        \draw[EdgeGraph](011)--(012);
        \draw[EdgeGraph](012)--(013);
        \draw[EdgeGraph](012)--(022);
        \draw[EdgeGraph](013)--(014);
        \draw[EdgeGraph](013)--(023);
        \draw[EdgeGraph](014)--(024);
        \draw[EdgeGraph](022)--(023);
        \draw[EdgeGraph](023)--(024);
    \end{tikzpicture}}
    \label{subfig:hill_poset_2_3}}
    \qquad \qquad
    \subfloat[][$\SetHill_\mathbf{1}(4)$]{
    \centering
    \scalebox{1}{
    \begin{tikzpicture}[Centering,xscale=.8,yscale=.8,
        x={(0,-1cm)},
        y={(-1cm,-1cm)},
        z={(1cm,-1cm)}]
        \foreach \x in {0,...,1} {
            \foreach \y in {0,...,2} {
                \foreach \z in {0,..., 3} {
                    \draw[LineGrid](\x,0,\z)--(\x,2,\z);
                    \draw[LineGrid](0,\y,\z)--(1,\y,\z);
                }
                \foreach \z in {0,..., 2} {
                    \draw[LineGrid](\x,\y,\z)--(\x,\y,\z+1);
                }
            }
        }
        \node[NodeGraph,MarkBA](0000)[]at(0,0,0){};
        \node[NodeGraph](0001)[]at(0,0,1){};
        \node[NodeGraph](0002)[]at(0,0,2){};
        \node[NodeGraph](0003)[]at(0,0,3){};
        \node[NodeGraph](0011)[]at(0,1,1){};
        \node[NodeGraph](0012)[]at(0,1,2){};
        \node[NodeGraph](0013)[]at(0,1,3){};
        \node[NodeGraph](0022)[]at(0,2,2){};
        \node[NodeGraph](0023)[]at(0,2,3){};
        \node[NodeGraph](0111)[]at(1,1,1){};
        \node[NodeGraph](0112)[]at(1,1,2){};
        \node[NodeGraph](0113)[]at(1,1,3){};
        \node[NodeGraph](0122)[]at(1,2,2){};
        \node[NodeGraph,MarkBA](0123)[]at(1,2,3){};
        \node[NodeLabeldGraph,left of=0000]{$0000$};
        \node[NodeLabeldGraph,left of=0001]{$0001$};
        \node[NodeLabeldGraph,left of=0002]{$0002$};
        \node[NodeLabeldGraph,left of=0003]{$0003$};
        \node[NodeLabeldGraph,left of=0011]{$0011$};
        \node[NodeLabeldGraph,left of=0012]{$0012$};
        \node[NodeLabeldGraph,left of=0013]{$0013$};
        \node[NodeLabeldGraph,left of=0022]{$0022$};
        \node[NodeLabeldGraph,left of=0023]{$0023$};
        \node[NodeLabeldGraph,left of=0111]{$0111$};
        \node[NodeLabeldGraph,left of=0112]{$0112$};
        \node[NodeLabeldGraph,left of=0113]{$0113$};
        \node[NodeLabeldGraph,left of=0122]{$0122$};
        \node[NodeLabeldGraph,left of=0123]{$0123$};
        \draw[EdgeGraph](0000)--(0001);
        \draw[EdgeGraph](0001)--(0002);
        \draw[EdgeGraph](0001)--(0011);
        \draw[EdgeGraph](0002)--(0003);
        \draw[EdgeGraph](0002)--(0012);
        \draw[EdgeGraph](0003)--(0013);
        \draw[EdgeGraph](0011)--(0012);
        \draw[EdgeGraph](0011)--(0111);
        \draw[EdgeGraph](0012)--(0013);
        \draw[EdgeGraph](0012)--(0022);
        \draw[EdgeGraph](0012)--(0112);
        \draw[EdgeGraph](0013)--(0023);
        \draw[EdgeGraph](0013)--(0113);
        \draw[EdgeGraph](0022)--(0023);
        \draw[EdgeGraph](0022)--(0122);
        \draw[EdgeGraph](0023)--(0123);
        \draw[EdgeGraph](0111)--(0112);
        \draw[EdgeGraph](0112)--(0113);
        \draw[EdgeGraph](0112)--(0122);
        \draw[EdgeGraph](0113)--(0123);
        \draw[EdgeGraph](0122)--(0123);
    \end{tikzpicture}}
    \label{subfig:hill_poset_1_4}}
    \caption{\footnotesize Hasse diagrams of some $\delta$-hill posets.}
    \label{fig:examples_hill_posets}
\end{figure}
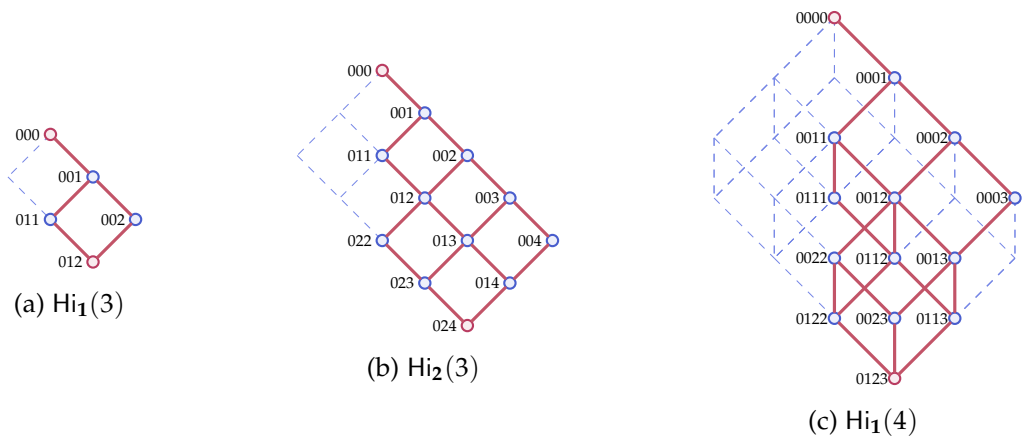

Let $\SetHill_\delta$ be the graded subset of $\SetCliff_\delta$ containing all
$\delta$-cliffs $u$ such that for any $i \in [|u| - 1]$, $u_i \leq u_{i + 1}$. Any element
of $\SetHill_\delta$ is a \Def{$\delta$-hill}.  Figure~\ref{fig:examples_hill_posets} shows
some ${\mathbf m}$-hill posets.
The ${\mathbf 1}$-hill posets are sometimes called Stanley lattices~\cite{Knu04}. The
$\delta$-hill posets can be seen as generalizations of these structures.

\begin{Proposition} \label{prop:properties_hill_objects}
    For any $m \geq 0$, $\SetHill_{\mathbf m}$ is straight, coated, closed by prefix,
    maximally extendable. Moreover, for any $n \geq 0$, $\SetHill_{\mathbf m}(n)$ is graded,
    a sublattice of $\SetCliff_{\mathbf m}(n)$, EL-shellable, nested, and constructible by
    interval doubling.
\end{Proposition}

\begin{Proposition} \label{prop:image_elevation_map_hill}
    For any range map $\delta$ and any $n \geq 0$,
    \begin{math}
        \ElevationImage_{\SetHill_\delta}(n)
        =
        \SetAvalanche_\delta(n).
    \end{math}
\end{Proposition}

A consequence of Propositions~\ref{prop:image_elevation_map_hill}
and~\ref{prop:enumeration_avalanche}, and of Point~\eqref{item:elevation_map_1} of
Proposition~\ref{prop:elevation_map} is that $\# \SetHill_{\mathbf m}(n) =
\FussCatalan_m(n)$.

\begin{Proposition} \label{prop:output_wings_hill}
    For any $m \geq 1$, $\OutputWings\Par{\SetHill_{\mathbf m}}$ is the set of the
    ${\mathbf m}$-cliffs $u$ satisfying
    \begin{math}
        u_1 \leq u_2 < \dots < u_{|u|},
    \end{math}
    and for all $i \in [2, |u|]$, $u_i < {\mathbf m}(i)$.
\end{Proposition}

\begin{Proposition} \label{prop:output_wings_hill_hill_bijection}
     For any $m \geq 1$ and $n \geq 0$, the map
    \begin{math}
        \phi : \OutputWings\Par{\SetHill_{\mathbf m}}(n)
        \to \SetHill_{\mathbf m - 1}(n)
    \end{math}
    defined, for any $u \in \OutputWings\Par{\SetHill_{\mathbf m}}(n)$ and $i \in [n]$, by
    \begin{math}
        \phi(u)_i := \IndicatorFunction_{i \ne 1} \Par{u_i - i + 2}
    \end{math}
    is an isomorphism of posets.
\end{Proposition}

\subsection{Canyon lattices} \label{sec:canyon_lattices}
\begin{figure}[ht]
    \centering
    \subfloat[][$\SetCanyon_\mathbf{1}(3)$]{
    \centering
    \scalebox{1}{
    \begin{tikzpicture}[Centering,xscale=.8,yscale=.8,rotate=-135]
        \draw[Grid](0,0)grid(1,2);
        \node[NodeGraph,MarkBA](000)[]at(0,0){};
        \node[NodeGraph](001)[]at(0,1){};
        \node[NodeGraph](002)[]at(0,2){};
        \node[NodeGraph](010)[]at(1,0){};
        \node[NodeGraph,MarkBA](012)[]at(1,2){};
        \node[NodeLabeldGraph,left of=000]{$000$};
        \node[NodeLabeldGraph,left of=001]{$001$};
        \node[NodeLabeldGraph,left of=002]{$002$};
        \node[NodeLabeldGraph,left of=010]{$010$};
        \node[NodeLabeldGraph,left of=012]{$012$};
        \draw[EdgeGraph](000)--(001);
        \draw[EdgeGraph](000)--(010);
        \draw[EdgeGraph](001)--(002);
        \draw[EdgeGraph](002)--(012);
        \draw[EdgeGraph](010)--(012);
    \end{tikzpicture}}
    \label{subfig:canyon_poset_1_3}}
    \qquad \qquad
    \subfloat[][$\SetCanyon_\mathbf{2}(3)$]{
    \centering
    \scalebox{1}{
    \begin{tikzpicture}[Centering,xscale=.8,yscale=.8,rotate=-135]
        \draw[Grid](0,0)grid(2,4);
        \node[NodeGraph,MarkBA](000)[]at(0,0,0){};
        \node[NodeGraph](001)[]at(0,1){};
        \node[NodeGraph](002)[]at(0,2){};
        \node[NodeGraph](003)[]at(0,3){};
        \node[NodeGraph](004)[]at(0,4){};
        \node[NodeGraph](010)[]at(1,0){};
        \node[NodeGraph](012)[]at(1,2){};
        \node[NodeGraph](013)[]at(1,3){};
        \node[NodeGraph](014)[]at(1,4){};
        \node[NodeGraph](020)[]at(2,0){};
        \node[NodeGraph](023)[]at(2,3){};
        \node[NodeGraph,MarkBA](024)[]at(2,4){};
        \node[NodeLabeldGraph,left of=000]{$000$};
        \node[NodeLabeldGraph,left of=001]{$001$};
        \node[NodeLabeldGraph,left of=002]{$002$};
        \node[NodeLabeldGraph,left of=003]{$003$};
        \node[NodeLabeldGraph,left of=004]{$004$};
        \node[NodeLabeldGraph,left of=010]{$010$};
        \node[NodeLabeldGraph,left of=012]{$012$};
        \node[NodeLabeldGraph,left of=013]{$013$};
        \node[NodeLabeldGraph,left of=014]{$014$};
        \node[NodeLabeldGraph,left of=020]{$020$};
        \node[NodeLabeldGraph,left of=023]{$023$};
        \node[NodeLabeldGraph,left of=024]{$024$};
        \draw[EdgeGraph](000)--(001);
        \draw[EdgeGraph](000)--(010);
        \draw[EdgeGraph](001)--(002);
        \draw[EdgeGraph](002)--(003);
        \draw[EdgeGraph](002)--(012);
        \draw[EdgeGraph](003)--(004);
        \draw[EdgeGraph](003)--(013);
        \draw[EdgeGraph](004)--(014);
        \draw[EdgeGraph](010)--(012);
        \draw[EdgeGraph](010)--(020);
        \draw[EdgeGraph](012)--(013);
        \draw[EdgeGraph](013)--(014);
        \draw[EdgeGraph](013)--(023);
        \draw[EdgeGraph](014)--(024);
        \draw[EdgeGraph](020)--(023);
        \draw[EdgeGraph](023)--(024);
    \end{tikzpicture}}
    \label{subfig:canyon_poset_2_3}}
    \qquad \qquad
    \subfloat[][$\SetCanyon_\mathbf{1}(4)$]{
    \centering
    \scalebox{1}{
    \begin{tikzpicture}[Centering,xscale=.8,yscale=.8,
        x={(0,-1cm)},
        y={(-1cm,-1cm)},
        z={(1cm,-1cm)}]
        \foreach \x in {0,...,1} {
            \foreach \y in {0,...,2} {
                \foreach \z in {0,..., 3} {
                    \draw[LineGrid](\x,0,\z)--(\x,2,\z);
                    \draw[LineGrid](0,\y,\z)--(1,\y,\z);
                }
                \foreach \z in {0,..., 2} {
                    \draw[LineGrid](\x,\y,\z)--(\x,\y,\z+1);
                }
            }
        }
        \node[NodeGraph,MarkBA](0000)[]at(0,0,0){};
        \node[NodeGraph](0001)[]at(0,0,1){};
        \node[NodeGraph](0002)[]at(0,0,2){};
        \node[NodeGraph](0003)[]at(0,0,3){};
        \node[NodeGraph](0010)[]at(0,1,0){};
        \node[NodeGraph](0012)[]at(0,1,2){};
        \node[NodeGraph](0013)[]at(0,1,3){};
        \node[NodeGraph](0020)[]at(0,2,0){};
        \node[NodeGraph](0023)[]at(0,2,3){};
        \node[NodeGraph](0100)[]at(1,0,0){};
        \node[NodeGraph](0101)[]at(1,0,1){};
        \node[NodeGraph](0103)[]at(1,0,3){};
        \node[NodeGraph](0120)[]at(1,2,0){};
        \node[NodeGraph,MarkBA](0123)[]at(1,2,3){};
        \node[NodeLabeldGraph,left of=0000]{$0000$};
        \node[NodeLabeldGraph,left of=0001]{$0001$};
        \node[NodeLabeldGraph,left of=0002]{$0002$};
        \node[NodeLabeldGraph,left of=0003]{$0003$};
        \node[NodeLabeldGraph,left of=0010]{$0010$};
        \node[NodeLabeldGraph,left of=0012]{$0012$};
        \node[NodeLabeldGraph,left of=0013]{$0013$};
        \node[NodeLabeldGraph,left of=0020]{$0020$};
        \node[NodeLabeldGraph,left of=0023]{$0023$};
        \node[NodeLabeldGraph,left of=0100]{$0100$};
        \node[NodeLabeldGraph,left of=0101]{$0101$};
        \node[NodeLabeldGraph,left of=0103]{$0103$};
        \node[NodeLabeldGraph,left of=0120]{$0120$};
        \node[NodeLabeldGraph,left of=0123]{$0123$};
        \draw[EdgeGraph](0000)--(0001);
        \draw[EdgeGraph](0000)--(0010);
        \draw[EdgeGraph](0000)--(0100);
        \draw[EdgeGraph](0001)--(0002);
        \draw[EdgeGraph](0001)--(0101);
        \draw[EdgeGraph](0002)--(0003);
        \draw[EdgeGraph](0002)--(0012);
        \draw[EdgeGraph](0003)--(0013);
        \draw[EdgeGraph](0003)--(0103);
        \draw[EdgeGraph](0010)--(0012);
        \draw[EdgeGraph](0010)--(0020);
        \draw[EdgeGraph](0012)--(0013);
        \draw[EdgeGraph](0013)--(0023);
        \draw[EdgeGraph](0020)--(0023);
        \draw[EdgeGraph](0020)--(0120);
        \draw[EdgeGraph](0023)--(0123);
        \draw[EdgeGraph](0100)--(0101);
        \draw[EdgeGraph](0100)--(0120);
        \draw[EdgeGraph](0101)--(0103);
        \draw[EdgeGraph](0103)--(0123);
        \draw[EdgeGraph](0120)--(0123);
    \end{tikzpicture}}
    \label{subfig:canyon_poset_1_4}}
    \caption{\footnotesize Hasse diagrams of some $\delta$-canyon posets.}
    \label{fig:examples_canyon_posets}
\end{figure}
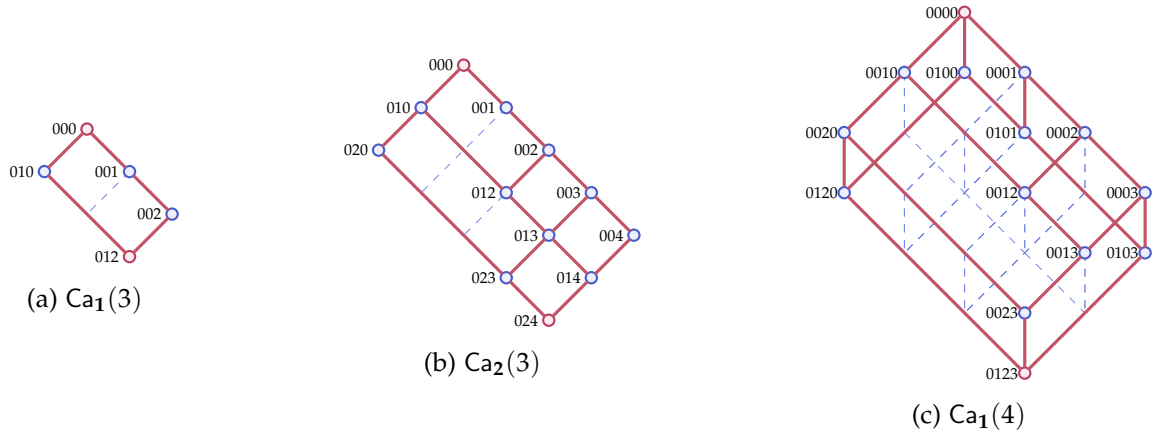

Let $\SetCanyon_\delta$ be the graded subset of $\SetCliff_\delta$ containing all
$\delta$-cliffs $u$ such that $u_{i - j} \leq u_i - j$, for all $i \in [|u|]$ and $j \in
\Han{u_i}$ satisfying $i - j \geq 1$.  Any element of $\SetCanyon_\delta$ is a
\Def{$\delta$-canyon}.  Figure~\ref{fig:examples_canyon_posets} shows some ${\mathbf
m}$-canyon posets.
The ${\mathbf 1}$-canyons are also know as Tamari diagrams and have been introduced
in~\cite{Pal86}. The set of these objects of size $n$ is in one-to-one correspondence with
the set of binary trees with $n$ internal nodes. It is also known that the componentwise
comparison of Tamari diagrams is the Tamari order~\cite{Pal86}. Moreover, for any $m \geq
2$, the ${\mathbf m}$-canyon posets are not the $m$-Tamari posets introduced
in~\cite{BPR12}. The $\delta$-canyon posets can be seen as different generalizations of
Tamari lattices.

\begin{Proposition} \label{prop:properties_canyon_objects}
    For any $m \geq 0$, $\SetCanyon_{\mathbf m}$ is straight, coated, closed by prefix, and
    minimally and maximally extendable. Moreover, for any $n \geq 0$, $\SetCanyon_{\mathbf
    m}(n)$ is a lattice, a meet semi-sublattice of $\SetCliff_{\mathbf m}(n)$, EL-shellable,
    nested, and constructible by interval doubling.
\end{Proposition}

Observe that $\SetCanyon_{\mathbf m}(n)$ is not a join semi-sublattice of $\SetCliff_{\mathbf
m}(n)$.  For instance,  the join operation of
$\SetCanyon_{\mathbf 1}(n)$ satisfies, by Theorem~\ref{thm:decr_incr_meet_join},
\begin{equation}
    00120 \JJoin_{\SetCanyon_{\mathbf 1}} 00201
    = \IncrMap_{\SetCanyon_{\mathbf 1}}(00120 \JJoin 00201)
    = \IncrMap_{\SetCanyon_{\mathbf 1}}(00221)
    = 00234.
\end{equation}
This computation of the join of two elements is similar to the ones described
in~\cite{Mar92} (see also~\cite{Gey94b}) for Tamari lattices.

\begin{Proposition} \label{prop:image_elevation_map_canyon}
    For any increasing range map $\delta$ and any $n \geq 0$,
    \begin{math}
        \ElevationImage_{\SetCanyon_\delta}(n)
        =
        \SetAvalanche_\delta(n).
    \end{math}
\end{Proposition}

\begin{Proposition} \label{prop:input_wings_butterflies_canyon}
    For any $m \geq 1$, $\InputWings\Par{\SetCanyon_{\mathbf m}}$ is the set of the
    ${\mathbf m}$-cliffs $u$ satisfying $u_i < u_{i + 1}$ for all $i \in [|u| - 1]$.
    Moreover, $\Butterflies\Par{\SetCanyon_{\mathbf m}}$ is the set of the ${\mathbf
    m}$-cliffs $u$ satisfying $1 \leq u_i < {\mathbf m}(i)$ for all $i \in [2, |u|]$, and
    $u_i - u_{i - 1} \geq 2$ for all $i \in [3, |u|]$.
\end{Proposition}

\begin{Theorem} \label{thm:poset_isomorphisms_wings_canyon_hill}
    For any $m \geq 1$ and $n \geq 0$,
    \begin{equation}
        \begin{tikzpicture}[Centering,xscale=4.6,yscale=1.4,font=\small]
            \node(ButterfliesCanyons)at(1,-1)
                {$\Butterflies\Par{\SetCanyon_{\mathbf m + 1}}(n)$};
            \node(InputWingCanyon)at(1,0){$\InputWings\Par{\SetCanyon_{\mathbf m}}(n)$};
            \node(Hill)at(2,0){$\SetHill_{\mathbf m - 1}(n)$};
            \node(OutputWingHill)at(2,-1){$\OutputWings\Par{\SetHill_{\mathbf m}}(n)$};
            \draw[Map](InputWingCanyon)--(ButterfliesCanyons)node[midway,left]{$\theta$};
            \draw[Map](InputWingCanyon)--(Hill)node[midway,above]{$\psi$};
            \draw[Map](OutputWingHill)--(Hill)node[midway,right]{$\phi$};
            \draw[Map](OutputWingHill)--(ButterfliesCanyons)node[midway,above]
                {$\theta \circ \psi^{-1} \circ \phi$};
        \end{tikzpicture}
    \end{equation}
    is a commutative diagram of isomorphisms of posets, where
    \begin{math}
        \theta : \InputWings\Par{\SetCanyon_{\mathbf m}}(n)
        \to \Butterflies\Par{\SetCanyon_{\mathbf m + 1}}(n)
    \end{math}
    and
    \begin{math}
        \psi : \InputWings\Par{\SetCanyon_{\mathbf m}}(n)
        \to \SetHill_{\mathbf m - 1}(n)
    \end{math}
    are the maps defined by
    \begin{math}
        \theta(u)_i := \IndicatorFunction_{i \ne 1} \Par{u_1 + i - 2}
    \end{math}
    and
    \begin{math}
        \psi(u)_i := u_i - i + 1,
    \end{math}
   for any $u \in \InputWings\Par{\SetCanyon_{\mathbf m}}(n)$ and
    $i \in [n]$.
\end{Theorem}

As a consequence of Theorem~\ref{thm:poset_isomorphisms_wings_canyon_hill}
and Proposition~\ref{prop:image_elevation_map_canyon},
\begin{equation}
    \# \Butterflies\Par{\SetCanyon_{\mathbf m + 1}}(n)
    = \# \InputWings\Par{\SetCanyon_{\mathbf m}}(n)
    = \# \SetCanyon_{\mathbf m - 1}(n)
    = \FussCatalan_{m - 1}(n).
\end{equation}

\begin{Theorem} \label{thm:morphism_canyon_hill}
    For any increasing range map $\delta$ and any $n \geq 0$, the map
    \begin{math}
        \ElevationMap_{\SetHill_\delta}^{-1}
        \circ
        \ElevationMap_{\SetCanyon_\delta}
        : \SetCanyon_\delta(n) \to \SetHill_\delta(n)
    \end{math}
    is both a bijection and a poset morphism.
\end{Theorem}

Even if
\begin{math}
    \ElevationMap_{\SetHill_\delta}^{-1}
    \circ
    \ElevationMap_{\SetCanyon_\delta}
    : \SetCanyon_\delta(n) \to \SetHill_\delta(n)
\end{math}
is a bijection, this map is not a poset isomorphism. Moreover, as a consequence of
Theorem~\ref{thm:morphism_canyon_hill}, for any $n \geq 0$, $\SetHill_\delta(n)$ is an order
extension of $\SetCanyon_\delta(n)$. To summarize, since by
Proposition~\ref{prop:elevation_map}, $\SetCanyon_\delta(n)$ is an order extension of
$\SetAvalanche_\delta(n)$, the three families of Fuss-Catalan posets fit into the chain
\begin{equation}
    \begin{tikzpicture}[Centering,xscale=3.4,yscale=2.6,font=\small]
        \node(Avalanche)at(0,0){$\SetAvalanche_\delta(n)$};
        \node(Canyon)at(1,0){$\SetCanyon_\delta(n)$};
        \node(Hill)at(2,0){$\SetHill_\delta(n)$};
        \draw[Map](Avalanche)--(Canyon)node[midway,above]
            {$\ElevationMap_{\SetCanyon_\delta}^{-1}$};
        \draw[Map](Canyon)--(Hill)node[midway,above]
            {$\ElevationMap_{\SetHill_\delta}^{-1}
        \circ
        \ElevationMap_{\SetCanyon_\delta}$};
        \draw[Map](Avalanche) edge[bend right=48] node[midway,above]
            {$\ElevationMap_{\SetHill_\delta}^{-1}$} (Hill);
    \end{tikzpicture}
\end{equation}
of posets for the order extension relation.

\section{Cliff associative algebras and quotients} \label{sec:algebras}

\subsection{Cliff associative algebras}
In the sequel, all the considered algebraic structures have a field $\K$ of characteristic
zero as ground field.  For any range map $\delta$, let $\SpaceCliff_\delta$ be the linear
span of all $\delta$-cliffs. This space is graded and decomposes as
\begin{equation}
    \SpaceCliff_\delta
    = \bigoplus_{n \geq 0} \SpaceCliff_\delta(n),
\end{equation}
where $\SpaceCliff_\delta(n)$, $n \geq 0$, is the linear span of all $\delta$-cliffs of size
$n$. By definition, the set $\Bra{\BasisF_u : u \in \SetCliff_\delta}$ is a basis of
$\SpaceCliff_\delta$, and we shall refer to it as the \Def{fundamental basis} or as the
\Def{$\BasisF$-basis}.  For any $n \geq 0$, the \Def{$\delta$-reduction map} is the map
\begin{math}
    \Reduction_\delta : \N^n \to \SetCliff_\delta(n)
\end{math}
defined for any word $u \in \N^n$ and any $i \in [n]$ by
\begin{math}
    \Reduction_\delta(u)_i := \min \Bra{u_i, \delta(i)}.
\end{math}
For instance, $\Reduction_{\mathbf 1}(212066) = 012045$ and $\Reduction_{\mathbf 2}(212066)
= 012066$.  Let
\begin{math}
    \Product : \SpaceCliff_\delta \otimes \SpaceCliff_\delta \to \SpaceCliff_\delta
\end{math}
be the binary product defined, for any $u, v \in \SetCliff_\delta$, by
\begin{equation} \label{equ:product_cliff_algebra}
    \BasisF_u \Product \BasisF_v =
    \sum_{\substack{
        v' \in \Reduction_\delta^{-1}(v) \\
        u v' \in \SetCliff_\delta
    }}
    \BasisF_{u v'}.
\end{equation}
For instance, in $\SpaceCliff_{\mathbf 2}$,

\begin{minipage}[t][0cm][b]{.6\textwidth} \small
\begin{equation}
    \BasisF_{00} \Product \BasisF_{011} =
    \BasisF_{00011} + \BasisF_{00111} + \BasisF_{00211} +
    \BasisF_{00311} + \BasisF_{00411},
\end{equation}
\end{minipage}
\begin{minipage}[t][0cm][b]{.32\textwidth} \small
\begin{equation}
    \BasisF_{0} \Product \BasisF_{0} = \BasisF_{00} + \BasisF_{01} + \BasisF_{02},
\end{equation}
\end{minipage}

\noindent and in $\SpaceCliff_{01312^\omega}$, we have

\begin{minipage}[t][0cm][b]{.6\textwidth} \small
\begin{equation}
    \BasisF_{00} \Product \BasisF_{011} =
    \BasisF_{00011} + \BasisF_{00111} + \BasisF_{00211}
    + \BasisF_{00311},
\end{equation}
\end{minipage}
\begin{minipage}[t][0cm][b]{.3\textwidth} \small
\begin{equation}
    \BasisF_{00} \Product \BasisF_{013} = 0.
\end{equation}
\end{minipage}

\begin{Theorem} \label{thm:cliff_associative_coalgebra}
    The space $\SpaceCliff_\delta$ endowed with the product $\Product$ is a magmatic unital
    graded algebra. Moreover, $\Product$ is associative if and only if $\delta$ is
    valley-free.
\end{Theorem}

We now establish a link between this product $\Product$ on the $\BasisF$-basis of
$\SpaceCliff_\delta$ and the $\delta$-cliff posets.  For this, let the two binary operations
$\Over$ and $\Under$ defined, for any $u, v \in \SetCliff_\delta$, by $u \Over v := u v$
and $u \Under v := u v'$ where $v'$ is the word on $\N$ of length $|v|$ satisfying, for any
$i \in [|v|]$,
\begin{equation}
    v'_i =
    \IndicatorFunction_{v_i = \delta(i)} \delta(|u| + i)
    + \IndicatorFunction_{v_i \ne \delta(i)} v_i.
\end{equation}
For instance, for $\delta = 112334^\omega$, $010 \Over 1021 = 010 1021$ and $010 \Under 1021
= 010 3041$.

\begin{Theorem} \label{thm:product_cliff_intervals}
    For any range map $\delta$, the product $\Product$ of $\SpaceCliff_\delta$ satisfies,
    for any $u, v \in \SetCliff_\delta$,
    \begin{equation}
        \BasisF_u \Product \BasisF_v =
        \IndicatorFunction_{u \Over v \in \SetCliff_\delta}
        \sum_{\substack{
            w \in \SetCliff_\delta \\
            u \Over v \Leq w \Leq u \Under v
        }}
        \BasisF_w.
    \end{equation}
\end{Theorem}

By mimicking the construction of bases of several combinatorial spaces by using a particular
partial order on their basis element (see for instance~\cite{DHT02,HNT05}), let for any $u
\in \SetCliff_\delta$,

\begin{minipage}[t][0cm][b]{.45\textwidth} \small
\begin{equation} \label{equ:basis_e_cliff}
    \BasisE_u :=
    \sum_{\substack{v \in \SetCliff_\delta \\
        u \Leq v
    }}
    \BasisF_v,
\end{equation}
\end{minipage}
\begin{minipage}[t][0cm][b]{.45\textwidth} \small
\begin{equation} \label{equ:basis_h_cliff}
    \BasisH_u :=
    \sum_{\substack{v \in \SetCliff_\delta \\
        v \Leq u
    }}
    \BasisF_v.
\end{equation}
\end{minipage}

\noindent By triangularity, the sets $\Bra{\BasisE_u : u \in \SetCliff_\delta}$ and
$\Bra{\BasisH_u : u \in \SetCliff_\delta}$ are bases of $\SpaceCliff_\delta$, called
respectively \Def{elementary basis} and \Def{homogeneous basis}, or respectively
\Def{$\BasisE$-basis} and \Def{$\BasisH$-basis}. For instance, in
$\SpaceCliff_{1021^\omega}$,

\begin{minipage}[t][0cm][b]{.42\textwidth} \small
\begin{equation}
    \BasisE_{1001} =
    \BasisF_{1001} + \BasisF_{1011} + \BasisF_{1021},
\end{equation}
\end{minipage}
\begin{minipage}[t][0cm][b]{.49\textwidth} \small
\begin{equation}
    \BasisH_{1001} =
    \BasisF_{1001} + \BasisF_{1000} + \BasisF_{0001} + \BasisF_{0000}.
\end{equation}
\end{minipage}

\begin{Proposition} \label{prop:product_cliff_e_h_basis}
    For any range map $\delta$, the product $\Product$ of $\SpaceCliff_\delta$ satisfies,
    for any $u, v \in \SetCliff_\delta$,

    \begin{minipage}[t][0cm][b]{.42\textwidth} \small
    \begin{equation}
        \BasisE_u \Product \BasisE_v =
        \IndicatorFunction_{u \Over v \in \SetCliff_\delta}
        \BasisE_{u \Over v},
    \end{equation}
    \end{minipage}
    \begin{minipage}[t][0cm][b]{.42\textwidth} \small
    \begin{equation}
        \BasisH_u \Product \BasisH_v
        = \BasisH_{\Reduction_\delta\Par{u \Under v}}.
    \end{equation}
    \end{minipage}
\end{Proposition}

A nonempty $\delta$-cliff $u$ is \Def{$\delta$-prime} if the relation $u = v \Over w$ with
$v, w \in \SetCliff_\delta$ implies $u = v$ or $u = w$. The graded collection of all these
elements is denoted by $\SetPrime_\delta$. For instance, for $\delta := 021^\omega$, among
others, the $\delta$-cliffs $0$, $01$, and $021111$ are $\delta$-prime, and $0210$ ($= 021
\Over 0$) is not.  Let $\AlphabetVar_{\SetPrime_\delta}$ be the alphabet $\Bra{a_u : u \in
\SetPrime_\delta}$ and $\K \Angle{\AlphabetVar_{\SetPrime_\delta}}$ be the free
(noncommutative) associative algebra generated by $\AlphabetVar_{\SetPrime_\delta}$. For any
$u \in \SetCliff_\delta$, we denote by $a^u$ the monomial $a_{u^{(1)}} \dots a_{u^{(k)}}$
where $u$ decomposes uniquely as $u = u^{(1)} \Over \dots \Over u^{(k)}$ where the
$u^{(i)}$, $i \in [k]$, are $\delta$-primes.  Let $\LeqSuffix$ be the partial order relation
on the monomials of $\K \Angle{\AlphabetVar_{\SetPrime_\delta}}$ wherein for any monomials
$a_{u^{(1)}} \dots a_{u^{(k)}}$ and $a_{v^{(1)}} \dots a_{v^{(\ell)}}$ of $\K
\Angle{\AlphabetVar_{\SetPrime_\delta}}$, one has
\begin{math}
    a_{u^{(1)}} \dots a_{u^{(k)}}
    \LeqSuffix a_{v^{(1)}} \dots a_{v^{(\ell)}}
\end{math}
if the word $u^{(1)} \dots u^{(k)}$ is a suffix of $v^{(1)} \dots v^{(\ell)}$. Given a set
$M$ of monomials of $\K \Angle{\AlphabetVar_{\SetPrime_\delta}}$, we denote by
$\min_{\LeqSuffix} M$ the set of all minimal elements of $M$ w.r.t. the order
relation~$\LeqSuffix$.

\begin{Theorem} \label{thm:cliff_algebra_presentation}
    For any valley-free range map $\delta$, the unital associative algebra
    $\SpaceCliff_\delta$ is isomorphic to
    \begin{math}
        \K \Angle{\AlphabetVar_{\SetPrime_\delta}}/
        _{\RelationSpace_\delta}
    \end{math}
    where $\RelationSpace_\delta$ is the associative algebra ideal of $\SpaceCliff_\delta$
    generated by the set
    \begin{equation}
        \min_{\LeqSuffix} \Bra{a^u a_v : u \in \SetCliff_\delta,
        v \in \SetPrime_\delta,
        \mbox{ and } u v \notin \SetCliff_\delta}.
    \end{equation}
\end{Theorem}

\begin{Proposition} \label{prop:cliff_algebra_finite_presentation}
    Let $\delta$ be a valley-free range map. As a unital associative algebra,
    $\SpaceCliff_\delta$ is
    \begin{inparaenum}[(i)]
        \item \label{item:cliff_algebra_finite_presentation_1}
        free if and only if $\delta$ is weakly increasing;
        \item \label{item:cliff_algebra_finite_presentation_2}
        admits a finite number of generators and a finite number of nontrivial relations
        between the generators if and only if $\delta$ is $1$-dominated.
    \end{inparaenum}
\end{Proposition}

A consequence of the freeness of $\SpaceCliff_{\mathbf 1}$ is that $\SpaceCliff_{\mathbf 1}$
is isomorphic as a unital associative algebra to $\FQSym$~\cite{MR95,DHT02}, an associative
algebra on the linear span of all permutations. This follows from the fact that $\FQSym$ is
also free as a unital associative algebra, and that its Hilbert series is the same as the
one of~$\SpaceCliff_{\mathbf 1}$. Moreover, in~\cite{NT14}, the authors construct some
associative algebras $^m \FQSym$ as generalizations of $\FQSym$ whose bases are indexed by
objects being generalizations of permutations. The algebras $\SpaceCliff_{\mathbf m}$, $m
\geq 0$, can therefore be seen as other generalizations of~$\FQSym$, not isomorphic to
$^m \FQSym$ when $m \geq 2$.

\subsection{Hill and canyon associative algebras}
Given a graded subset $\SubFamilly$ of $\SetCliff_\delta$, let
\begin{math}
    \SpaceCliff_{\SubFamilly}
    :=
    \SpaceCliff_\delta /_{\SpaceV_\SubFamilly}
\end{math}
be the quotient space of $\SpaceCliff_\delta$ where $\SpaceV_\SubFamilly$ is the linear span
of the set
\begin{math}
    \Bra{\BasisF_u : u \in \SetCliff_\delta \setminus \SubFamilly}.
\end{math}
By definition,
\begin{math}
    \Bra{\BasisF_u : u \in \SubFamilly}
\end{math}
is a basis of $\SpaceCliff_\SubFamilly$.  The graded subset $\SubFamilly$ is \Def{closed by
suffix reduction} if for any $n \geq 0$ and any $u \in \SubFamilly (n)$, for all suffixes
$u'$ of $u$, $\Reduction_\delta\Par{u'} \in \SubFamilly$.

\begin{Proposition} \label{prop:quotient_cliff}
    Let $\delta$ be a valley-free range map and $\SubFamilly$ be a graded subset of
    $\SetCliff_\delta$. If $\SubFamilly$ is closed by prefix and is closed by suffix
    reduction, then $\SpaceCliff_\SubFamilly$ is a quotient algebra of the unital
    associative algebra~$\SpaceCliff_\delta$.
\end{Proposition}

We also say that $\SpaceCliff_\SubFamilly$ has the \Def{interval condition} if the support
of any product $\BasisF_u \Product \BasisF_v$, $u, v \in \SubFamilly$, is empty or is an
interval of a poset $\SubFamilly(n)$, $n \geq 0$. When for any $n \geq 0$, $\SubFamilly(n)$
is a lattice, we denote by $\Meet'$ (resp. $\JJoin'$) its meet (resp. join) operation.
In this case, $\SubFamilly$ is \Def{meet-stable} (resp. \Def{join-stable}) if,
for any $n \geq 0$ and any $u, v \in \SubFamilly(n)$, 
the relation $u_i = v_i$ for an $i \in [n]$ implies that the $i$-th
letter of $u \Meet' v$ (resp. $u \JJoin' v$) is equal to~$u_i$.

\begin{Theorem} \label{thm:quotient_cliff_bounds_f_basis}
    Let $\delta$ be a valley-free range map and $\SubFamilly$ be a graded subset of
    $\SetCliff_\delta$ closed by prefix and by suffix reduction. If at least one of the
    following conditions is satisfied:
    \begin{inparaenum}[(i)]
        \item \label{item:quotient_cliff_bounds_f_basis_1} for any $n \geq 0$, all posets
        $\SubFamilly(n)$ are sublattices of $\SetCliff_\delta(n)$;
        \item \label{item:quotient_cliff_bounds_f_basis_2} for any $n \geq 0$, all posets
        $\SubFamilly(n)$ are meet semi-sublattices of $\SetCliff_\delta(n)$, maximally
        extendable, and join-stable;
    \end{inparaenum}
    then $\SpaceCliff_\SubFamilly$ has the interval condition.
\end{Theorem}

\paragraph{Hill associative algebras.}
For any $m \geq 0$, let $\SpaceHill_m$ be the quotient $\SpaceCliff_{\SetHill_{\mathbf m}}$.
This quotient is well-defined due to the fact that $\SetHill_{\mathbf m}$ satisfies the
conditions of Proposition~\ref{prop:quotient_cliff}. Moreover, by
Proposition~\ref{prop:properties_hill_objects} and
Point~\eqref{item:quotient_cliff_bounds_f_basis_1} of
Theorem~\ref{thm:quotient_cliff_bounds_f_basis}, $\SpaceHill_m$ has the interval condition.
For instance, one has in $\SpaceHill_1$,

\begin{minipage}[t][0cm][b]{.6\textwidth} \small
\begin{equation}
    \BasisF_{01} \Product \BasisF_{01}
    =
    \BasisF_{0111} + \BasisF_{0112} + \BasisF_{0113} + \BasisF_{0122} + \BasisF_{0123},
\end{equation}
\end{minipage}
\begin{minipage}[t][0cm][b]{.3\textwidth} \small
\begin{equation}
    \BasisF_{01} \Product \BasisF_{00}
    =
    0.
\end{equation}
\end{minipage}

\noindent In $\SpaceHill_2$, one has

\begin{minipage}[t][0cm][b]{.52\textwidth} \small
\begin{equation}
    \BasisF_{02} \Product \BasisF_{023}
    =
    \BasisF_{02223} + \BasisF_{02233} + \BasisF_{02333},
\end{equation}
\end{minipage}
\begin{minipage}[t][0cm][b]{.35\textwidth} \small
\begin{equation}
    \BasisF_{011} \Product \BasisF_{01}
    =
    \BasisF_{01111}.
\end{equation}
\end{minipage}

\noindent For any $m \geq 1$, $\SpaceHill_m$ is not free as unital associative algebra.
Indeed, the quasi-inverse of the generating series of its generators (see
Table~\ref{tab:number_generators_cliff_quotients} for the first coefficients obtained with
the help of the computer) is not the Hilbert series of $\SpaceHill_m$, which is expected
when this algebra is free.
\begin{table}
    \centering
    \subfloat[][In $\SpaceHill_m$.]{
    \begin{scriptsize}
    \begin{tabular}{c|ccccccccc}
        $m$ & \multicolumn{9}{c}{Number of generators of $\SpaceHill_{\mathbf m}$} \\ \hline
        $1$ & 0 & 1 & 1 & 2 & 6  & 18  & 59  & 196 & 669 \\
        $2$ & 0 & 1 & 2 & 7 & 33 & 168 & 900 & 4980 & \\
    \end{tabular}
    \end{scriptsize}}
    \label{subtab:number_generators_cliff_quotients_hill}
    \qquad
    \subfloat[][In $\SpaceCanyon_m$.]{
    \begin{scriptsize}
    \begin{tabular}{c|ccccccccc}
        $m$ & \multicolumn{9}{c}{Number of generators of $\SpaceCanyon_{\mathbf m}$}
            \\ \hline
        $1$ & 0 & 1 & 1 & 2 & 5  & 14  & 42  & 132 & 429 \\
        $2$ & 0 & 1 & 2 & 7 & 30 & 149 & 788 & 4332 & \\
    \end{tabular}
    \end{scriptsize}}
    \label{subtab:number_generators_cliff_quotients_canyon}
    \caption{\footnotesize The first numbers, dimension by dimension, of generators in some
    quotients of~$\SpaceCliff_{\mathbf m}$.}
    \label{tab:number_generators_cliff_quotients}
\end{table}

\paragraph{Canyon associative algebras.}
For any $m \geq 0$, let $\SpaceCanyon_m$ be the quotient $\SpaceCliff_{\SetCanyon_{\mathbf
m}}$. This quotient is well-defined due to the fact that $\SetCanyon_{\mathbf m}$ satisfies
the conditions of Proposition~\ref{prop:quotient_cliff}. Moreover, by
Proposition~\ref{prop:properties_canyon_objects}, the fact that for any $m \geq 0$ and $n
\geq 0$, $\SetCanyon_{\mathbf m}(n)$ is join-stable, and by
Point~\eqref{item:quotient_cliff_bounds_f_basis_2} of
Theorem~\ref{thm:quotient_cliff_bounds_f_basis}, $\SpaceCanyon_m$ has the interval
condition. For instance, one has in $\SpaceCanyon_1$,

\begin{minipage}[t][0cm][b]{.4\textwidth} \small
\begin{equation}
    \BasisF_{0} \Product \BasisF_{01}
    =
    \BasisF_{001} + \BasisF_{002} + \BasisF_{012},
\end{equation}
\end{minipage}
\begin{minipage}[t][0cm][b]{.45\textwidth} \small
\begin{equation}
    \BasisF_{0} \Product \BasisF_{002}
    =
    \BasisF_{0002} + \BasisF_{0003} + \BasisF_{0103}.
\end{equation}
\end{minipage}

\noindent In $\SpaceCanyon_2$, one has

\begin{minipage}[t][0cm][b]{.45\textwidth} \small
\begin{equation}
    \BasisF_{01} \Product \BasisF_{0014}
    =
    0,
\end{equation}
\end{minipage}
\begin{minipage}[t][0cm][b]{.45\textwidth} \small
\begin{equation}
    \BasisF_{01} \Product \BasisF_{0013}
    =
    \BasisF_{010013}.
\end{equation}
\end{minipage}

\noindent The associative algebra $\SpaceCanyon_1$ is the Loday-Ronco algebra~\cite{LR98},
also known as $\PBT$~\cite{HNT05}. Moreover, for any $m \geq 2$, $\SpaceCanyon_m$ is
not free as unital associative algebra.  Indeed, the quasi-inverse of the generating series
of its generators (see Table~\ref{tab:number_generators_cliff_quotients} for the first
coefficients) is not the Hilbert series of $\SpaceCanyon_m$, which is expected when this
algebra is free.

\begin{footnotesize}
\bibliographystyle{plain}
\bibliography{Bibliography}
\end{footnotesize}

\end{document}